\documentclass[12pt]{article}
\usepackage{amsmath}
\usepackage{amssymb}
\usepackage{mathtools}
\usepackage{dsfont}
\usepackage{cite}
\newsymbol \blackbox 1004
\newcommand{\eh}{\hfill}\newlength{\sperr}

\newenvironment{proof}{{\settowidth{\sperr}{\bf\rm
Proof}%
\par\addvspace{0.3cm}\noindent\parbox[t]{1.3\sperr}
{\bf\rm P\eh r\eh o\eh o\eh f\eh }%
}}{\nopagebreak\mbox{}
$\blackbox$\par\addvspace{0.3cm}}

\def\nn{\nonumber}

\def\a{\alpha}
\def\b{\beta}

\def\G{\Gamma}

\def\la{\lambda}
\def\om{\omega}

\def\th{\theta}

\def\vp{\varphi}
\def\vt{\vartheta}

\def\wh{\widehat}
\def\wt{\widetilde}
\def\ov{\overline}

\def\BC{{\mathbb C}}

\def\cla{{\mathcal A}}

\def\cll{{\mathcal L}}

\def\clt{\mathcal{T}}

\def\mfa{\mathfrak{A}}

\def\diag{\mathrm{diag}}
\def\col{\mathrm{col}}

\newcommand{\I}{\mathrm{i}}
\newcommand{\1}{\mathbf{1}}

\newtheorem{Pa}{Paper}[section]
\newtheorem{Tm}[Pa]{{\bf Theorem}}

\newtheorem{Cy}[Pa]{{\bf Corollary}}
\newtheorem{Rk}[Pa]{{\bf Remark}}

\newtheorem{Pn}[Pa]{{\bf Proposition}}

\title{On  the inversion of the block double-structured and of the   triple-structured Toeplitz
matrices and on the corresponding reflection coefficients}

\author{Inna Roitberg and Alexander Sakhnovich}

\date{}
\begin{document}
\maketitle

\begin{abstract} 
The results on the inversion of convolution operators as well as Toeplitz (and block Toeplitz) matrices
in the $1$-D (one-dimensional) case are classical and have numerous
applications. Last year, we considered the $2$-D case of Toeplitz-block Toeplitz (TBT) matrices,
described a minimal information, which is necessary to recover the inverse
matrices, and gave a complete characterisation of the inverse matrices.
Now, we develop our approach for the more complicated cases of block TBT-matrices and $3$-D
Toeplitz matrices. 
\end{abstract}

\vspace{0.5em}


\vspace{0.5em}


\section{Introduction}
\setcounter{equation}{0}
 
The well-known Toeplitz matrices $T$ are diagonal-constant matrices, that is, they have the form 
\begin{align}\label{R1}&
T=\{\clt_{i-k}\}_{i,k=1}^n, \quad \clt_r\in \BC,
\end{align}
where $\BC$ stands for the complex plain. The  theory of Toeplitz matrices and operators (and structured operators in general) is an important part of analysis closely
related to many other mathematical domains and to various applications (see, e.g., \cite{BoG, BoeS, GS, Kac, Pol, SaL, Sar, Sar1, Sim} and references
therein). The inversion of Toeplitz matrices and their continuous analogs (i.e., convolution operators) 
is an essential part of this theory. Note that  structured operators appear as a result of certain  homogeneities in the studied processes.
For instance, Toeplitz matrices and convolution operators play crucial role in stationary, stable and Levy processes
in statistics (see, e.g., \cite{SaL15}). Multi-structured operators appear in the study of the processes depending on several variables.

The inversion of convolution operators on the real line, semi-axis and finite intervals is connected with the names of N.~Wiener, E. Hopf, N. Levinson, M.G.~Krein,
 I.C. Gohberg, V.A. Ambartsumyan,  L.A. Sakhnovich and many other mathematicians and applied scientists 
 (see the results and references in \cite{Amb,  BrH, Dev, Go, GoKr,  Kre, Levin, Nag, Nob, SaSaR, SaL15, Shi, Tr, Wid}).  
 For the closely related and actively studied  inversion of the Toeplitz matrices see, for instance, \cite{GoS, HR, SaA73, SaA80, SaSp, KKM0, KKM, KS, Bin0}
 and references therein, including important works by T. Kailath and coauthors and by D. Bini and coauthors.
 A fruitful method of operator identities (see \cite{SaL68,  SaL73, SaL86, SaL, SaL15} and references therein) was applied to the inversion of 
 Toeplitz and block Toeplitz matrices and study of their properties in  \cite{SaA73, SaA80, SaAJFA, SaSp}. In particular, 
 the structure of the inverses to Toepltz matrices was derived in \cite{SaA73}.
 Further developments in the inversion of the so called structured matrices and various fast methods of inversion are connected
with the  works \cite{Bin0, HR, KKM0, KKM, KS, KaS2}. Interesting matrix identities (or, equivalently, Lyapunov and Sylvester  equations) appeared
in those works. For instance, the equations 
$T-ZTZ^*=Q, \quad Z:=\{\delta_{i, k+1}\}_{i,k=1}^n$
were considered in the  papers
\cite{KKM0, KKM}.


When the entries $\clt_r$ in \eqref{R1} are Toeplitz blocks
(instead of being scalars) one talks about Toeplitz-block Toeplitz  matrices $T$. {\it Later, we use for  Toeplitz-block Toeplitz  matrices the acronym TBT.} 
One may consider TBT-matrices as the $2$-D (two-dimensional) analog of the Toeplitz matrices. 
In spite of a number of important recent and older works \cite{DaJ, CheB, GeK,   Jus, Kab, Kalo, Lev, Nap, Olsh, WaK} on the inversion of the TBT-matrices 
and  of the convolution operators in multidimensional spaces, the structure of the corresponding inverse matrices and operators  was first characterised
in our recent  paper  \cite{ALS-Trans} on  the TBT-matrix case. In particular, the fact that   some important space-time covariance matrices have the Toeplitz-block Toeplitz form
stimulated the study of the  inverted  multi-dimensional Toeplitz matrices. Among various other  interesting recent works on  the  Toeplitz matrices, convolution
operators and their applications, we mention \cite{AGKLS, Bin, BogoB,  CIL, DerSi, GeS, JeVa, KreL, Pest0, Pest, Xi}.

In this paper, we consider an important (and more complicated than TBT-matrices) case of block TBT-matrices:
\begin{align}\label{D0}&
T=\{\clt_{i-k}\}_{i,k=1}^{m_1}, \quad \clt_r=\{t^{(r)}_{i-k}\}_{i,k=1}^{m_2}, 
\end{align}
where $t^{(r)}_s$ are $m_3 \times m_3$ blocks (not necessarily Toeplitz) and  $m_p \geq 2$. A special
subclass of the block TBT-matrices such that  $\quad U_3 t_j^{(r)}U_3=(t_j^{(r)})^{\tau}$ $(U_3:=\{\delta_{m_3-i-k+1}\}_{i,k=1}^{m_3})$,
where $t^{\tau}$ is   the  transpose of the block $t$,
{\it is denoted by the acronym DSTU} (double structured Toeplitz matrices with the property ``U").  DSTU includes, in particular, $3$-D Toeplitz matrices:
\begin{align}\label{D1}&
T=\{\clt_{i-k}\}_{i,k=1}^{m_1}, \quad \clt_r=\{t^{(r)}_{i-k}\}_{i,k=1}^{m_2}, \quad t^{(r)}_s=\{\tau^{(r,s)}_{i-k}\}_{i,k=1}^{m_3}, \quad \tau^{(r,s)}_{i-k}\in \BC.
\end{align}
The subclasses of the self-adjoint block TBT-matrices and of the DSTU-matrices  are of essential interest and admit more complete
analysis than general block TBT-matrices.

Clearly, matrices $T$ of the forms \eqref{D0} and \eqref{D1} are $m\times m$ matrices, where
\begin{align}\label{D1'}&
m:=m_1m_2m_3 \quad (m_p\geq 2).
\end{align}
It is easy to see that a $3$-D Toeplitz matrix $T$ satisfies three matrix identities of the form 
\begin{align}\label{D4}&
A_pT-TA_p^*=\I\big(M_{1p}M_{2p}+ M_{3p}M_{4p}\big),
 \end{align}
where $A_p^*$ is the conjugate transpose of the matrix $A_p$,
\begin{align}
\label{D3'}&
A_1=\cla_1 \otimes I_{m_2m_3}=\cla_1 \otimes I_{m_2}\otimes I_{m_3},
\\ \label{D3}&
 A_2=I_{m_1} \otimes( \cla_2\otimes I_{m_3})=(I_{m_1} \otimes \cla_2)\otimes I_{m_3}, \quad A_3=I_{m_1m_2}\otimes\cla_3;
\\ \label{D2}&
\cla_p= \{a_{j-\ell}\}_{j,\ell=1}^{m_p}\quad (p=1,2,3), 
\quad a_r=\left\{\begin{array}{l}0 \,\, {\mathrm{for}}\,\,  r<0, \\
\I / 2  \,\, {\mathrm{for}}\,\,  r=0, \\ \I , \,\, {\mathrm{for}}\,\,  r>0; \end{array} \right. 
 \end{align}
$I_r$ is the $r\times r$ identity matrix, $\otimes$ stands for the Kronecker product;
$M_{1p}$ and $M_{3p}$ are $m \times \frac{m}{m_p}$ matrices, and 
$M_{2p}$ and $M_{4p}$ are $\frac{m}{m_p} \times m$ matrices,
see the definitions of $M_{kp}$ in \eqref{D13}--\eqref{D12} and \eqref{D5}--\eqref{D8}.
We note that the second index in $M_{kp}$ shows the number of the corresponding
matrix identity and the first index shows the place of the matrix on the right-hand side of this
identity (differently from, e.g., $G_{kp}$ which are the blocks of $G$).
{\it The block TBT matrices \eqref{D0} satisfy two matrix identities, namely, identities \eqref{D4}, where $p=1,\, 2$.}

It easily follows from  the mixed-product property of the Kronecker products that the matrices $A_p$ $(p=1,2,3)$ given by \eqref{D3'} and \eqref{D3} commute pairwise.
In view of the mixed-product property (and the bilinearity), we also have:
\begin{align}
\label{D3!}&
(A_1-zI)^{-1}=(\cla_1-I_{m_1})^{-1} \otimes I_{m_2}\otimes I_{m_3},
 \\
\label{D3++}&
 (A_2-zI)^{-1}=I_{m_1} \otimes( \cla_2-zI_{m_2})^{-1}\otimes I_{m_3}, \\
\label{D3+}&
(A_3-zI)^{-1}=I_{m_1}\otimes I_{m_2}\otimes (\cla_3-zI_{m_3})^{-1}.
 \end{align}
Further in the text, we consecutively use Kronecker products and usually do not mention their standard properties
 including bilinearity, associativity,  mixed-product property
and equalities 
$$F(I_k\otimes h)=F\otimes h, \quad F(h\otimes I_k)=h\otimes F,$$ 
as well as the adjoint equalities 
$$(I_k\otimes h^*)F^*=F^*\otimes h^*, \quad (h^*\otimes I_k)F^*=h^*\otimes F^*,$$
 where the matrices $F$ have $k$ columns and $h$ are
row vectors. The resolvents of $\cla_p$ are known explicitly and further we will need the following equality (see, for instance, \cite[(1.10)]{SaAJFA}):
\begin{align}\label{F0}&
\begin{bmatrix} 1 & 1 & \ldots & 1\end{bmatrix}\big(\cla_p^*-z I_{m_p}\big)^{-1}
=-\frac{2}{2z+\I}\begin{bmatrix} 1 & \frac{2z-\I}{2z+\I} & \ldots & \left(\frac{2z-\I}{2z+\I}\right)^{m_p-1}\end{bmatrix}.
\end{align}

When the block TBT-matrix $T$ is invertible, relations \eqref{D4} yield
\begin{align}\label{I1}&
RA_p-A_p^*R=\I R\Pi_p \wh \Pi_p R,
 \end{align}
where $p=1,2$ and
\begin{align}\label{I2}&
 R:=T^{-1}, \quad  \Pi_p:=\begin{bmatrix}M_{1p} & M_{3p}\end{bmatrix}, \quad
\wh \Pi_p:=\begin{bmatrix}M_{2p} \\ M_{4p}\end{bmatrix}.
\end{align}
The inverse matrix $R$ is explicitly recovered from each of the identities \eqref{I1} $(p=1,2)$ if one has ``information" about $R\Pi_p$ and  $\wh \Pi_p R$.
However, the block TBT-matrix is determined by  its $(2m_1-1)(2m_2-1)(m_3)^2$ entries and the matrices $R\Pi_p$ and  $\wh \Pi_p R$
have together $4m^2/m_p$ entries. In other words, too much ``information" is required in order to recover $R$ from one identity.
{\it Here, the main task and the main difficulty is to use both identities \eqref{I1} simultaneously and  to minimise  the required
``information" as well as to get in this way a better understanding of the structure of $R$.}

 For this purpose, together with $T^{-1}$ we   consider  the so-called matrix reflection coefficient
 \begin{align}\label{D21}&
\rho(x,y)=h({y})^{\, \tau}T^{-1}h(x),
 \end{align}
where $x=(x_1,x_2)$, $y=(y_1,y_2)$,  and $h$ is expressed via column vectors $h_p$:
\begin{align}\label{D22}&
h(x):=h_1(x_1)\otimes h_2(x_2) \otimes I_{m_3}, \quad h_p(x_p):=\{x_p^{i-1})\}_{i=1}^{m_p}.
\end{align}
Recall that $h^{\, \tau}$  is the  transpose of $h$.
{\it The function   $\rho(x,y)$ is an $m_3 \times m_3$ matrix polynomial
which uniquely determines
$T^{-1}$ $($in a simple way$)$ and is of interest in itself.}
Using \eqref{D3!}, \eqref{D3++} and \eqref{F0}, we show  that 
$\rho(x,y)$ is easily expressed via
\begin{align}\label{D24}&
\om(\la,\mu)=\cll^*\Big(\prod_{p=1}^2(A_p^*-\mu_p I)^{-1}\Big)T^{-1}\Big(\prod_{p=1}^2(A_p-\la_p I)^{-1}\Big)\cll,
 \end{align}
where  $\la=(\la_1,\la_2)$, $\mu=(\mu_1,\mu_2)$, 
\begin{align}\label{D25}&
 \cll=\1_{m_1m_2}\otimes I_{m_3}, \quad \1_{r}:=\col\begin{bmatrix}1 & 1 &\ldots & 1\end{bmatrix},
 \end{align}
col means column and $\1_r$ is an $r\times 1$ column vector. Indeed, we have
$$\cll^*\Big(\prod_{p=1}^2(A_p^*-\mu_p I)^{-1}\Big)=\1_{m_1}^*(\cla_{1}-\mu_1I_{m_1})^{-1}\otimes \1_{m_2}^*(\cla_{2}-\mu_2I_{m_2})^{-1}\otimes I_{m_3}.$$
Thus, the equality 
\begin{align}\label{I3}&
\cll^*\Big(\prod_{p=1}^2(A_p^*-\mu_p I)^{-1}\Big)= -(y_1-1)(y_2-1)h({y})^{\, \tau}
\end{align}
holds for $y_p=\frac{\mu_p-\I/2}{\mu_p+\I/2}$ ($p=1,2$) or, equivalently, for
\begin{align}\label{I4}&
\mu_p=-\vp(y_p), \quad \vp(z)=\frac{\I}{2}\frac{z+1}{z-1}.
\end{align}
Similar considerations are valid for $\Big(\prod_{p=1}^2(A_p-\la_p I)^{-1}\Big)\cll$, and so we derive
\begin{equation}\label{I5}
\rho(x,y)=\prod_{p=1}^2(x-x_p)^{-1}\prod_{j=1}^2(y-y_j)^{-1}\om\big(\vp(x_1),\vp(x_2),-\vp(y_1),-\vp(y_2)\big).
\end{equation}
Hence, the study of $T^{-1}$ is reduced to the study of $\om(\la,\mu)$.

Some notations were already introduced in the introduction.
The notation diag stands for the diagonal or block diagonal matrix (e.g. $\diag\{d_1,\ldots, d_n\} $ stands
for the diagonal matrix with $d_1,\ldots, d_n$ on the main diagonal).


\section{Block TBT-matrices: preliminaries}\label{Prel}
\setcounter{equation}{0}
Let $T$ be an invertible block TBT matrix \eqref{D0}. Similarly to \cite[(2.11)--(2.17)]{ALS-Trans}, we obtain
the expressions for the matrices $M_{kp}$ $(p=1,2)$ such that the matrix identities \eqref{D4} with $p=1,2$ hold.
Namely, we have
 \begin{align}\label{D13}&
M_{11}=\col  \begin{bmatrix} M_{11}^{(1)} & M_{11}^{(2)} & \ldots & M_{11}^{(m_1)}\end{bmatrix}, \quad 
M_{11}^{(i)}:=\frac{1}{2}\clt_0 +\sum_{s=1}^{i-1}\clt_{s},
\\ \label{D14}&
M_{21}= \1_{m_1}^*\otimes I_{m_2m_3}=1_{m_1}^*\otimes I_{m_2}\otimes I_{m_3}, \quad M_{31}=M_{21}^*,
\\ \label{D15}&
M_{41}=  \begin{bmatrix} M_{41}^{(1)} & M_{41}^{(2)} & \ldots & M_{41}^{(m_1)}\end{bmatrix}, \quad 
M_{41}^{(k)}:=\frac{1}{2}\clt_0 +\sum_{s=1}^{k-1}\clt_{-s};
\end{align}
and
\begin{align}\label{D9}&
M_{12}=\{M_{12}^{(i-k)}\}_{i,k=1}^{m_1}, \quad M_{22}=I_{m_1} \otimes \1_{m_2}^*\otimes I_{m_3},
\\ \label{D10}&
  M_{32}=M_{22}^*, \quad 
M_{42}=\{M_{42}^{(i-k)}\}_{i,k=1}^{m_1},
\\ \label{D11}&
 M_{12}^{(r)}=\col \begin{bmatrix} \frac{1}{2}t_0^{(r)} &  ({t_0^{(r)}}/{2})+t_1^{(r)}&\ldots & ({t_0^{(r)}}/{2})+\displaystyle{\sum_{j=1}^{m_2-1}}t_{j}^{(r)}\end{bmatrix},
\\ \label{D12}&
M_{42}^{(r)}=\begin{bmatrix} \frac{1}{2}t_0^{(r)} & ({t_0^{(r)}}/{2})+t_{-1}^{(r)}& \ldots & ({t_0^{(r)}}/{2})+\displaystyle{\sum_{\ell=1}^{m_2-1}}t^{(r)}_{-\ell}\end{bmatrix};
\end{align}
where $t_i^{(r)}$ are $m_3\times m_3$ blocks.

Further  in the text, we usually assume that $T$ is  invertible, and we set
\begin{align}\label{D17}&
R:=T^{-1}; \quad  \G_p:=R\Pi_p, \quad \wh \G_p:=\wh \Pi_p R;
\\ \label{D17+}& \Pi_p:=\begin{bmatrix}M_{1p} & M_{3p}\end{bmatrix}, \quad
\wh \Pi_p:=\begin{bmatrix}M_{2p} \\ M_{4p}\end{bmatrix};
\\ \label{D18}&
\G:=\begin{bmatrix} \G_1 & \G_2\end{bmatrix}, \quad 
 \wh \G:=\col \begin{bmatrix} \wh \G_1 & \wh \G_2\end{bmatrix},
 \end{align}
 where $p=1,2$ (and some of the notations were already shortly discussed in  the introduction).
 Identities \eqref{D4} may be rewritten in the form
  \begin{align}\label{D16}&
A_pT-TA_p^*=\I \Pi_p \wh \Pi_p, 
 \end{align}
which yields  
   \begin{align}\label{D19}&
R A_p-A_p^*R=\I \G_p \wh \G_p.
 \end{align}
 \begin{Pn} \label{simprep} Let $\om(\la, \mu)$ be defined by \eqref{D24}, \eqref{D25}. Then, we have
  \begin{align}\label{D26}&
\om(\la, \mu)= \I (\la_p-\mu_p)^{-1} u_p(\mu)  \wh u_p(\la) \quad (p=1,2),
 \end{align}
 where
   \begin{align} \label{D28} &
  u_1(\mu)=v_1(\mu)-\I\big(\1_{m_2}^*(\cla_2^*-\mu_2I_{m_2})^{-1}\otimes I_{m_3}\big)\begin{bmatrix} I_{m_2m_3} & 0\end{bmatrix},
  \\ &  \label{D28+}
 \wh   u_1(\la)=\wh v_1(\la) +\I \begin{bmatrix}0 \\ I_{m_2m_3} \end{bmatrix}\big((\cla_2-\la_2I_{m_2})^{-1}\1_{m_2}\otimes I_{m_3}\big);
   \end{align}
 and
 \begin{align}
   \label{D27} &
  u_2(\mu)=v_2(\mu)-\I\big(\1_{m_1}^*(\cla_1^*-\mu_1I_{m_1})^{-1}\otimes I_{m_3}\big)\begin{bmatrix} I_{m_1m_3} & 0\end{bmatrix},
  \\ &  \label{D27+}
 \wh   u_2(\la)=\wh v_2(\la) +\I \begin{bmatrix}0 \\ I_{m_1m_3} \end{bmatrix}\big((\cla_1-\la_1I_{m_1})^{-1}\1_{m_1}\otimes I_{m_3}\big),
\\ & \label{D28++}
v_p(\mu):=\cll^*\Big(\prod_{i=1}^2(A_i^*-\mu_i I)^{-1}\Big)\G_p, \quad \wh v_p(\la):=\wh \G_p\Big(\prod_{i=1}^2(A_i-\la_i I)^{-1}\Big)\cll.
 \end{align} 
 \end{Pn}
\begin{proof}. In view of \eqref{D19}, we have 
 \begin{align}\label{D29-}&
(\mu_p-\la_p)R= R(A_p-\la_pI)-(A_p^*-\mu_p I)R-\I \G_p \wh \G_p,
 \end{align}
or, equivalently,
 \begin{align}\nn
 (A_p^*-\mu_p I)^{-1}R(A_p-\la_pI)^{-1}= &(\mu_p-\la_p)^{-1}\big((A_p^*-\mu_p I)^{-1}R-R(A_p-\la_pI)^{-1}
\\ \label{D29} & 
-\I (A_p^*-\mu_p I)^{-1}\G_p \wh \G_p(A_p-\la_pI)^{-1}\big).
 \end{align}
 Next, we set 
  \begin{align}\label{D33} & 
\cll_p:=\1_{m_p}\otimes I_{m_3},
 \end{align}
and easily obtain (for $\cll$ given in \eqref{D25}) the equalities
   \begin{align}\label{D34} & 
\cll_2^*M_{21}=\cll_1^*M_{22}=\cll^*.
 \end{align}
 In view of \eqref{D3!}, \eqref{D3++}  and \eqref{D14}, \eqref{D9}, using again  the mixed-product property of the Kronecker product we derive
   \begin{align}\label{D35} & 
M_{21}(A_1^*-zI)^{-1}=\1_{m_1}^*(\cla_1^*-zI_{m_1})^{-1}\otimes I_{m_2}\otimes I_{m_3}, 
\\ \label{D36} & 
 M_{22}(A_2^*-zI)^{-1}=I_{m_1}
\otimes \1_{m_2}^*(\cla_2^*-zI_{m_2})^{-1}\otimes I_{m_3}.
 \end{align}
 Moreover, we have representations
   \begin{align}\label{D37} & 
M_{21}=\begin{bmatrix} I_{m_2m_3} & \ldots & I_{m_2m_3}\end{bmatrix}, \quad M_{22}=\diag\{\cll_2^*, \cll_2^*,\ldots, \cll_2^*\}.
 \end{align}
It follows from \eqref{D33}--\eqref{D35} and from \eqref{D37} that 
   \begin{align}\nn
\cll^* (A_1^*-zI)^{-1}&=\cll_2^*\Big(\big(\1_{m_1}^*(\cla_1^*-zI_{m_1})^{-1}\otimes I_{m_2}\big)\otimes I_{m_3}\Big)
\\ & \nn
=\Big(1_{m_2}^*\big(\1_{m_1}^*(\cla_1^*-zI_{m_1})^{-1}\otimes I_{m_2}\big)\Big)\otimes I_{m_3}
\\ & \nn
=\1_{m_1}^*(\cla_1^*-zI_{m_1})^{-1}\otimes \1_{m_2}^*\otimes I_{m_3}=\1_{m_1}^*(\cla_1^*-zI_{m_1})^{-1}\otimes \cll_2^*
\\ & \label{D37'} 
=\big(\1_{m_1}^*(\cla_1^*-zI_{m_1})^{-1}\otimes I_{m_3}\big)M_{22}.
 \end{align} 
 In a similar way, relations \eqref{D33}, \eqref{D34} and \eqref{D36}, \eqref{D37} yield
 \begin{align}\nn
\cll^* (A_2^*-zI)^{-1}&=\cll_1^*\Big(\big(I_{m_1}\otimes \1_{m_2}^*(\cla_2^*-zI_{m_2})^{-1}\big)\otimes I_{m_3}\Big)
\\ & \nn
=\Big(\1_{m_1}^*\big(I_{m_1}\otimes \1_{m_2}^*(\cla_2^*-zI_{m_2})^{-1}\big)\Big)\otimes I_{m_3}
\\ & \nn
=\1_{m_1}^*\otimes \1_{m_2}^*(\cla_2^*-zI_{m_2})^{-1}\otimes I_{m_3}
\\ & \nn
=\big( \1_{m_2}^*(\cla_2^*-zI_{m_2})^{-1}\otimes I_{m_3}\big)\big(\1_{m_1}^*\otimes I_{m_2m_3}\big)
\\ & \label{D38} 
=\big( \1_{m_2}^*(\cla_2^*-zI_{m_2})^{-1}\otimes I_{m_3}\big)M_{21}.
 \end{align} 
 Setting in \eqref{D29} $p=2$, multiplying both parts of \eqref{D29} by $\cll^*  (A_1^*-\mu_1 I)^{-1}$ from the left
 and by $ (A_1^*-\la_1 I)^{-1}\cll$ from the right, and using \eqref{D37'}, we obtain
   \begin{align}\nn
\om(\la, \mu)= &\I (\la_2-\mu_2)^{-1} \Big(\cll^*\Big(\prod_{i=1}^2(A_i^*-\mu_i I_m)^{-1}\Big)\G_2 \wh \G_2\Big(\prod_{i=1}^2(A_i-\la_i I_m)^{-1}\Big)\cll
\\ & \nn
-\I\big(\1_{m_1}^*(\cla_1^*-\mu_1I_{m_1})^{-1}\otimes I_{m_3}\big)\begin{bmatrix} I_{m_1m_3} & 0\end{bmatrix}\wh \G_2\Big(\prod_{i=1}^2(A_i-\la_i I_m)^{-1}\Big)\cll
\\ & \nn
+\I \cll^*\Big(\prod_{i=1}^2(A_i^*-\mu_i I_m)^{-1}\Big)\G_2\begin{bmatrix}0 \\ I_{m_1m_3} \end{bmatrix}\big((\cla_1-\la_1I_{m_1})^{-1}\1_{m_1}\otimes I_{m_3}\big)\Big).
 \end{align}
 Here, we took into account the equality $M_{32}=M_{22}^*$ and the definitions \eqref{D17}, \eqref{D17+}.
 The formula for $\om$ above shows that \eqref{D26} holds for $p=2$ and for $u_2, \, \wh u_2$ given by
 \eqref{D27}, \eqref{D27+} and \eqref{D28++}. In the same way (using \eqref{D29} and \eqref{D38}) one can show that \eqref{D26} is valid for $p=1$.
 \end{proof}
{\it Proposition \ref{simprep} is an analog of the celebrated $1$-D result for convolution operators} \cite{Amb, Sob} (see also \cite{SaL15} for the corresponding $1$-D formula
in a general case and further discussion). Next, we generalise the main result from \cite{ALS-Trans}.
\section{Block TBT-matrices: main results}\label{Block}
\setcounter{equation}{0}
{\bf 1.}  The basic  role in our considerations is played by the  square \\ $2(m_1+m_2)m_3\times 2(m_1+m_2)m_3$ matrix function $G(\la)$:
 \begin{align}\label{g}&
G(\la):=\begin{bmatrix}G_{11}(\la_2) & G_{12} \\
G_{21} & G_{22}(\la_1)
\end{bmatrix},
\\ \label{g1}&
\sbox0{$\begin{matrix}(\cla_1-\la_1 I_{m_1})\otimes I_{m_3} & 0\\ 0 & (\cla_1-\la_1 I_{m_1})\otimes I_{m_3}\end{matrix}$}
G_{22}(\la_1)=\left[\usebox{0} \right],
\\ \label{g2}&\sbox1{$\begin{matrix} (\cla_2-\la_2 I_{m_2})\otimes I_{m_3} & 0\\ 0 & (\cla_2-\la_2 I_{m_2})\otimes I_{m_3}\end{matrix}$}
G_{11}(\la_2)=\left[\usebox{1} \right],
\\ \label{g3}&
G_{21} =\I\left(\wh \Pi_2 \G_1-\begin{bmatrix} \cll_1\cll_2^* & 0 \\ K_{11} & 0\end{bmatrix}\right), \quad
G_{12} =\I\left(\wh \Pi_1 \G_2-\begin{bmatrix} \cll_2\cll_1^* & 0 \\ K_{12} & 0\end{bmatrix}\right),
 \end{align} 
where $\wh \Pi_i \G_p=\wh \Pi_i  T^{-1} \Pi_p$; $G_{21}$ and $G_{12}$, respectively,  are  $2m_1m_3\times 2m_2m_3$ and $2m_2m_3\times 2m_1m_3$
constant matrices; the matrices $\cll_p$ are given by \eqref{D33}, and $K_{1p}$ are given by
\begin{align} 
 \label{g5}&
K_{11}=\col \begin{bmatrix}\frac{1}{2}M_{42}^{(0)} & \frac{1}{2}M_{42}^{(0)}+M_{42}^{(1)} & \ldots & \frac{1}{2} M_{42}^{(0)}+\sum_{j=1}^{m_1-1}M_{42}^{(j)}
\end{bmatrix},
\\ \label{g4}&
K_{12}= \begin{bmatrix}\frac{1}{2}M_{12}^{(0)} & \frac{1}{2}M_{12}^{(0)}+M_{12}^{(-1)} & \ldots & \frac{1}{2} M_{12}^{(0)}+\sum_{j=1}^{m_1-1}M_{12}^{(-j)}
\end{bmatrix}.
\end{align}
We note that $G_{22}$ depends on $\la_1$ (and $\cla_1$) and 
$G_{11}$ depends on $\la_2$ (and $\cla_2$).
The matrices $K_{1p}$ appear in the matrix identities \eqref{D45} for $M_{4p}$. 

We set
\begin{align}\label{V1}&
 u(\mu) =\begin{bmatrix}u_1(\mu) & u_2(\mu)
\end{bmatrix}, \quad \wh u(\la) = \begin{bmatrix}\wh u_1(\la) \\ \wh u_2(\la)
\end{bmatrix}.
\end{align}
\begin{Tm}\label{MbTm} Let the block TBT-matrix $T$ be invertible. Then, the matrix function $\wh u(\la)$ of the form \eqref{V1}, 
where $\wh u_p$ $(p=1,2)$ are defined in Proposition~\ref{simprep} and appear in the formulas \eqref{D26} for $\om(\la,\mu)$, 
is given by the relations
\begin{align}\label{T1}&
\wh u(\la) = G(\la)^{-1}\col \begin{bmatrix}0 & \1_{m_2}\otimes I_{m_3}  & 0 &  \1_{m_1}\otimes I_{m_3}
\end{bmatrix}\th(\la),
\\  \label{T2}&
\th(\la):=\I\Big(I_{m_3}+KR   \Big(\prod_{i=1}^2(A_i-\la_i I)^{-1}\Big)\cll              \Big),
\\  \label{T3}&
K:=  \begin{bmatrix}\frac{1}{2}M_{42}^{(0)} & \frac{1}{2}M_{42}^{(0)}+M_{42}^{(-1)} &\ldots & \frac{1}{2} M_{42}^{(0)}+\sum_{j=1}^{m_1-1}M_{42}^{(-j)}
\end{bmatrix}.
\end{align}
\end{Tm}
\begin{Rk}\label{RkTau} The matrix function $G(\la)$ is determined by its constant blocks  $G_{12}$ and $G_{21}$. In the case of TBT-matrices,
$G_{21}$ is easily expressed via $G_{12}$ and vice versa. The function $u$ in \eqref{D26} is easily expressed
via $\wh u$ in that case as well. These interrelations between $G_{21}$ and $G_{12}$ and between $u$ and $\wh u$ are caused by the 
interrelations of the \eqref{F2} type between $T$ and its transpose $T^{\tau}$. Such interrelations do not hold in the general block TBT-case. The situation is in some
respects similar to the
differences between {\rm 1-{D}} Toeplitz and block Toeplitz cases. In particular, one cannot effectively simplify
the equalities
$$(\la_2-\mu_2)u_1(\mu)\wh u_1(\mu)=(\la_1-\mu_1)u_2(\mu)\wh u_2(\mu),
$$
implied by \eqref{D26}.
However,  simple interrelations between $G_{21}$ and $G_{12}$ and between $u$ and $\wh u$ 
appear again, when one considers the self-adjoint or DSTU subclasses.
\end{Rk}
\begin{Rk}\label{RkTheta}
Note that although $\th(\la)$  is expressed in terms of $G$ in the TBT-case \cite{ALS-Trans},
the situation in the block TBT-case is more complicated even for the self-adjoint or DSTU subclasses.
\end{Rk}
{\bf Open Problem I.} Find the cases of block TBT-matrices, where $\th(\la)$ is uniquely determined by $G(\la)$. 
\begin{proof}  of Theorem \ref{MbTm}.
In addition to the formulas \eqref{D35} and \eqref{D36} for the matrices $M_{21}$ and $M_{22}$, taking into account  \eqref{D3!} and \eqref{D3++})
we derive the following relations, which we write (for convenience) for the adjoint matrices $M_{31}$ and $M_{32}$:
\begin{align}\nn
M_{31}\big((\cla_2-\la_2I_{m_2})^{-1}\otimes I_{m_3}\big)&=\1_{m_1}\otimes (\cla_2-\la_2I_{m_2})^{-1}\otimes I_{m_3}
\\ \label{D39}&
=(A_2-\la_2 I)^{-1}M_{31},
\\ \label{D40}
M_{32}\big((\cla_1-\la_1I_{m_1})^{-1}\otimes I_{m_3}\big)&=(A_1-\la_1 I)^{-1}M_{32}.
\end{align} 
Next,  we substitute $\mu_p=\la_p$ into \eqref{D29-} and (in view of the obtained identity and of the definitions of $\G_p$, $\wh u_p$ and $\cll_p$) we have
\begin{align}\label{D41}
(A_p^*-\la_pI)^{-1}\G_p\wh u_p(\la)=&(A_p^*-\la_p I)^{-1}\G_p\wh \G_p\Big(\prod_{i=1}^2(A_i-\la_i I_m)^{-1}\Big)\cll
\\ & \nn
+\I(A_p^*-\la_p I)^{-1}RM_{3p}\big((\cla_k-\la_k I_{m_k})^{-1}\1_{m_k}\otimes I_{m_3}\big)
\\ = & \nn
\I R\Big(\prod_{i=1}^2(A_i-\la_i I_m)^{-1}\Big)\cll
\\ & \nn
- \I (A_p^*-\la_p I)^{-1}R
(A_k-\la_kI)^{-1}\cll
\\ & \nn
+\I(A_p^*-\la_p I)^{-1}RM_{3p}\big((\cla_k-\la_k I_{m_k})^{-1}\otimes I_{m_3}\big)\cll_k,
\end{align}
where $k,p\in \{1,2\}$, $k\not=p$ (i.e., {\it the equalities hold for $p=1, \,\, k=2$ and for  $p=2, \,\, k=1$}).

 Using the equalities $M_{31}\cll_2=M_{32}\cll_1=\cll$ (which are immediate
from \eqref{D34}) and relations \eqref{D39}, \eqref{D40}, we simplify the right-hand side of  \eqref{D41}:
\begin{align}\label{D42}&
(A_p^*-\la_pI)^{-1}\G_p\wh u_p(\la)=\I R\Big(\prod_{i=1}^2(A_i-\la_i I_m)^{-1}\Big)\cll .
\end{align}
Introduce the matrix functions $\mfa_p(\la_p)$ by the equalities
\begin{align}\label{mfa}&
\mfa_p(\la_p)=\begin{bmatrix}(\cla_p^*-\la_p I_{m_p})\otimes I_{m_3} & 0\\ 0 & (\cla_p-\la_p I_{m_p})\otimes I_{m_3}\end{bmatrix} \quad (p=1,2),
 \end{align}
and the matrix function $F(\la)$ by the equality
\begin{equation}\label{D44}
F(\la)=\begin{bmatrix} \mathfrak{A}_2(\la_2) &\I \mathfrak{A}_2(\la_2) \wh \Pi_1 (A_2^*-\la_2 I)^{-1}\G_2
\\
 \I \mathfrak{A}_1(\la_1) \wh \Pi_2 (A_1^*-\la_1 I)^{-1}\G_1
  &\mathfrak{A}_1(\la_1)
 \end{bmatrix}.
 \end{equation} 
 In view of the second equalities in \eqref{D28++} and \eqref{V1}, and relations \eqref{D42} and \eqref{D44}, we have
 $$F(\la)\wh u(\la)=\begin{bmatrix}   \mfa_2(\la_2)\big(-\wh v_1(\la)+\wh u_1(\la)\big)\\ \mfa_1(\la_1)\big(\wh u_2(\la) -\wh v_2(\la)\big)
 \end{bmatrix},$$
which, taking into account \eqref{D28+}, \eqref{D27+} and \eqref{mfa}, reduces to the equality 
\begin{align}\label{D43}&
F(\la)\wh u(\la) =\I \, \col \begin{bmatrix}0 & \1_{m_2}\otimes I_{m_3}  & 0 &  \1_{m_1}\otimes I_{m_3}
\end{bmatrix}.
\end{align}

Besides commutation properties \eqref{D39}, \eqref{D40} for the matrices $M_{3p}$, we will need matrix identities
for $M_{4p}$. Indeed, according to \eqref{D0}, \eqref{D15} and \eqref{D9} (for $ k,p \in \{1,2\}, \,\, k\not=p$) we have 
\begin{align}\label{D45}&
(\cla_k\otimes I_{m_3})M_{4p}-M_{4p}A_k^*=\I Q_k=\I (K_{1k}K_{2k}+K_{3k}K_{4k}),
\end{align}
where
\begin{align} \label{D46}&
K_{11}=\col \begin{bmatrix}\frac{1}{2}M_{42}^{(0)} & \frac{1}{2}M_{42}^{(0)}+M_{42}^{(1)} & \ldots & \frac{1}{2} M_{42}^{(0)}+\sum_{j=1}^{m_1-1}M_{42}^{(j)}
\end{bmatrix},
\\ \label{D47}& 
K_{21}=\1_{m_1}^*\otimes I_{m_2m_3}=M_{21}, \quad K_{31}=\1_{m_1}\otimes I_{m_3}, \quad K_{41}=K,
\\ \label{D48}&
K=  \begin{bmatrix}\frac{1}{2}M_{42}^{(0)} & \frac{1}{2}M_{42}^{(0)}+M_{42}^{(-1)} &\ldots & \frac{1}{2} M_{42}^{(0)}+\sum_{j=1}^{m_1-1}M_{42}^{(-j)}
\end{bmatrix};
\\ \label{D49}&
K_{12}= \begin{bmatrix}\frac{1}{2}M_{12}^{(0)} & \frac{1}{2}M_{12}^{(0)}+M_{12}^{(-1)} & \ldots & \frac{1}{2} M_{12}^{(0)}+\sum_{j=1}^{m_1-1}M_{12}^{(-j)}
\end{bmatrix},
\\ \label{D50}& 
K_{22}=I_{m_1}\otimes 1_{m_2}^*\otimes I_{m_3}=M_{22}, \quad K_{32}=\1_{m_2}\otimes I_{m_3}, \quad K_{42}=K.
\end{align}
Recall that $K_{11}$ and $K_{12}$ were already introduced in \eqref{g5} and \eqref{g4}, and $K$ was given in \eqref{T3}.
It easily follows from \eqref{D45}, \eqref{D47} and \eqref{D50} that
\begin{align}\label{D50+}&
Q_p=K_{1p}M_{2p}+\cll_p K.
\end{align}
Taking adjoints of the both parts in the equalities \eqref{D39} and \eqref{D40}, we obtain
\begin{align}\label{D51}&
M_{2p}(A_k^*-\la_k I)^{-1}=\big((\cla_k^*-\la_k I_{m_k})^{-1}\otimes I_{m_3}\big)M_{2p}.
\end{align}
Relations \eqref{mfa}, \eqref{D45} and \eqref{D51}, and the block representation of $\wh \Pi$ in \eqref{D17+} yield
an important equality
\begin{align}\label{D52}&
\wh \Pi_p (A_k^*-\la_k I)^{-1}=\mathfrak{A}_k(\la_k)^{-1}\left(\wh \Pi_p+\I
\begin{bmatrix}0 \\ Q_k\end{bmatrix} (A_k^*-\la_k I)^{-1}\right)
 \end{align}
for $ k,p \in \{1,2\}, \,\, k\not=p$. Clearly, the indices $k$ and $p$ may change places above. Next, multiplying \eqref{D42} by $Q_p$ we have
\begin{align}\label{D53}&
Q_p(A_p^*-\la_pI)^{-1}\G_p\wh u_p(\la)=\I Q_p R\Big(\prod_{i=1}^2(A_i-\la_i I)^{-1}\Big)\cll .
\end{align}
Using the definitions \eqref{D28+} and \eqref{D27+} of $\wh u_p$ and formula \eqref{D50+}, we rewrite the right-hand side of \eqref{D53}:
\begin{align}\label{D54}&
\I Q_p R\Big(\prod_{i=1}^2(A_i-\la_i I_m)^{-1}\Big)\cll =\I \begin{bmatrix} K_{1p} & 0\end{bmatrix}\wh u_p(\la)+\I \cll_p \a(\la),
\\ \label{D55}&
\a(\la):=KR\Big(\prod_{i=1}^2(A_i-\la_i I)^{-1}\Big)\cll.
\end{align}
For  $ k,p \in \{1,2\}, \,\, k\not=p$, formulas \eqref{D52}--\eqref{D54} yield
\begin{align} \nn
\wh \Pi_k (A_p^*-\la_p I)^{-1}\G_p\wh u_p(\la)=&\mathfrak{A}_p(\la_p)^{-1}\left(\wh \Pi_k\G_p-
\begin{bmatrix}0 & 0\\ K_{1p} & 0\end{bmatrix} \right)\wh u_p(\la)
\\ \label{D56}&
-\mathfrak{A}_p(\la_p)^{-1}\begin{bmatrix}0 \\ \cll_p\end{bmatrix} \a(\la).
 \end{align}  
 Taking into account \eqref{D56}, we rewrite \eqref{D44} and \eqref{D43} in the form
\begin{align}\label{D57}&
\wt F(\la)\wh u(\la) = \col \begin{bmatrix}0 & \cll_2  & 0 &  \cll_1
\end{bmatrix}\th(\la);
\\ \label{D58}& \th(\la)=\I\big(I_{m_3}+\a(\la)\big), \quad
\wt F(\la)=\begin{bmatrix} \mathfrak{A}_2(\la_2) & \wt F_{12}, 
\\
\wt F_{21}  &\mathfrak{A}_1(\la_1)
\end{bmatrix}, 
\\ 
 \label{D59}& 
\wt F_{21} = =\I\left(\wh \Pi_2 \G_1-\begin{bmatrix} 0 & 0 \\ K_{11} & 0\end{bmatrix}\right), \quad
\wt F_{12} =\I\left(\wh \Pi_1 \G_2-\begin{bmatrix} 0 & 0 \\ K_{12} & 0\end{bmatrix}\right),
\end{align}
where  $\cll_p$  are  given by \eqref{D33}.
In view of \eqref{g}--\eqref{g3} and \eqref{D58}, \eqref{D59}, we have
\begin{align}\label{D60}&
G(\la)=\wt F(\la)+\I \wh \cll \wh \cll^*, \quad \wh \cll:=\col \begin{bmatrix} \cll_2 & 0 & -\cll_1 & 0 \end{bmatrix}.
\end{align}
It follows from the definitions \eqref{D17}--\eqref{D18} and equalities \eqref{D34} that
\begin{align}\label{D61}&
\wh \cll^*\wh \G=(\cll_2^*M_{21}-\cll_1^*M_{22})R=0.
\end{align}
Now, the definitions \eqref{D28+}, \eqref{D27+}, \eqref{D28++}, and \eqref{V1} of $\wh u_p$ and $\wh u$ together with the equalities \eqref{D61} yield
\begin{align}\label{D62}&
\wh \cll^*\wh u(\la)=0.
\end{align} 
Using \eqref{D60} and \eqref{D62}, we rewrite \eqref{D57} in the form
\begin{align}\label{D63}&
G(\la)\wh u(\la) = \col \begin{bmatrix}0 & \cll_2  & 0 &  \cll_1
\end{bmatrix}\th(\la),
\end{align}
where $\th(\la)$ coincides with  $\th(\la)$ in \eqref{T2}. Hence, \eqref{T1} is immediate.
\end{proof}

{\bf 2.} In order to construct $u(\mu)$, we introduce matrix function
 \begin{align}\label{E}&
E(\mu):=\begin{bmatrix}E_{11}(\mu_2) & E_{12} \\
E_{21} & E_{22}(\mu_1)
\end{bmatrix},
\\ \label{E1}&
\sbox0{$\begin{matrix}(\cla_1^*-\mu_1 I_{m_1})\otimes I_{m_3} & 0\\ 0 & (\cla_1^*-\mu_1 I_{m_1})\otimes I_{m_3}\end{matrix}$}
E_{22}(\mu_1)=\left[\usebox{0} \right],
\\ \label{E2}&\sbox1{$\begin{matrix} (\cla_2^*-\mu_2 I_{m_2})\otimes I_{m_3} & 0\\ 0 & (\cla_2^*-\mu_2 I_{m_2})\otimes I_{m_3}\end{matrix}$}
E_{11}(\mu_2)=\left[\usebox{1} \right],
\\ \label{E3}&
E_{21} =-\I\left(\wh \Pi_2 \G_1-\begin{bmatrix} 0 & 0 \\ K_{11} & \cll_1\cll_2^*\end{bmatrix}\right), \quad
E_{12} =-\I\left(\wh \Pi_1 \G_2-\begin{bmatrix} 0 & 0 \\ K_{12} & \cll_2\cll_1^*\end{bmatrix}\right).
 \end{align} 
\begin{Tm} \label{bM2Tm} Let the block TBT-matrix $T$ be invertible. Then, the matrix function $u(\mu) =\begin{bmatrix}u_1(\mu) & u_2(\mu)
\end{bmatrix}$
where $u_p$ $(p=1,2)$ are defined in Proposition~\ref{simprep} and appear in the formulas \eqref{D26} for $\om(\la,\mu)$, 
is given by the relations
 \begin{align}\label{C13}&
 u(\mu) =\vt(\mu)\begin{bmatrix} \cll_2^* & 0 & \cll_1^* & 0\end{bmatrix}E(\mu)^{-1},
 \\ \label{C14}&
 \vt(\mu):=-\I\Big(I_{m_3}+ \cll^*\Big(\prod_{i=1}^2(A_i^*-\mu_i I_m)^{-1}\Big)RN\Big),
\\ \label{C5}&
N:=\col \begin{bmatrix} \frac{1}{2} M_{12}^{(0)} &   \frac{1}{2}M_{12}^{(0)}+ M_{12}^{(1)} &\ldots & 
\frac{1}{2}M_{12}^{(0)}+ \sum_{j=1}^{m_1-1}M_{12}^{(j)}\end{bmatrix}. 
 \end{align}   
\end{Tm}
\begin{proof}. Similar to the proof of \eqref{D42}, we take into account \eqref{D29-} (with $\mu_p=\la_p$), \eqref{D37'}, \eqref{D38} and, using  this time the
 definitions of $\wh \G_p$ and $u_p$ (instead of the definitions of $\G_p$ and $\wh u_p$), derive
\begin{align}\nn
u_p\wh \G_p (A_p-\mu_p I)^{-1}=&-\I \cll^*\Big(\prod_{i=1}^2(A_i^*-\mu_i I_m)^{-1}\Big)R
\\ &\nn
+\I \cll^*(A_k^*-\mu_k I_m)^{-1}R(A_p-\mu_p I)^{-1}
\\ & \nn
-\I \big(\1_{m_k}^*(\cla_{m_k}-\mu_k I_{m_k})^{-1}\otimes I_{m_3}\big)M_{2p}R(A_p-\mu_p I)^{-1}
\\ \label{C1} =&
-\I \cll^*\Big(\prod_{i=1}^2(A_i^*-\mu_i I_m)^{-1}\Big)R,
\end{align}
where $k,p\in \{1,2\}$, $k\not=p$. From \eqref{C1} we obtain an analog of the equality \eqref{D43}, namely, we have
\begin{align}\nn  &
u(\mu)\begin{bmatrix} \mathfrak{A}_2(\mu_2) &-\I \wh \G_1 (A_1-\mu_1 I)^{-1}\Pi_2 \mathfrak{A}_1(\mu_1)
\\
 -\I \wh \G_2 (A_2-\mu_2 I)^{-1}\Pi_1 \mathfrak{A}_2(\mu_2) 
   &\mathfrak{A}_1(\mu_1)
 \end{bmatrix}
 \\ & \label{C2}
 =-\I \begin{bmatrix}\1_{m_2}^*\otimes I_{m_3} & 0  &  \1_{m_1}^*\otimes I_{m_3}  &  0
\end{bmatrix}.
 \end{align} 
 In view of \eqref{D39} and \eqref{D40}, an analog of \eqref{D52} takes the form
 \begin{align}\label{C3}&
 (A_p-\mu_p I)^{-1}\Pi_k=\left( \Pi_k-\I (A_p-\mu_p I)^{-1}
\begin{bmatrix} V_p& 0\end{bmatrix} \right)\mathfrak{A}_p(\mu_p)^{-1},
 \end{align}
 where $\I V_p$ is the right-hand side of the identity
 \begin{align}\label{C3+}&
A_pM_{1k}-M_{1k}\big(\cla_p^*\otimes I_{m_3}\big)=\I V_p.
 \end{align} 
 Similar to \eqref{D50+}, we derive
 \begin{align}\label{C4}&
 V_p=N\cll_p^*+M_{3p}K_{1k},
 \end{align}  
 where $k,p\in \{1,2\}$ and $k\not=p$; $\cll_p$ are defined in \eqref{D34}, $M_{3p}$ are given in \eqref{D14} and \eqref{D10}, $K_{1k}$ are introduced in \eqref{g5} and \eqref{g4},
 and $N$ has the form \eqref{C5}. 

 Next, in view of \eqref{D28++}, \eqref{C1} and \eqref{C4}, we note that
  \begin{align}\nn
u_p (\mu)\wh \G_p(A_p-\mu_p I)^{-1} V_p &=-\I \cll^*\Big(\prod_{i=1}^2(A_i^*-\mu_i I_m)^{-1}\Big)R\Pi_p\begin{bmatrix}0\\ K_{1k}\end{bmatrix}-\I \b(\mu)\cll_p^*
\\ \label{C6}&
=-\I u_p(\mu)\begin{bmatrix}0\\ K_{1k}\end{bmatrix}-\I \b(\mu)\cll_p^*,
\end{align}  
where
\begin{align}\label{C7}&
 \b(\mu)= \cll^*\Big(\prod_{i=1}^2(A_i^*-\mu_i I_m)^{-1}\Big)RN.
 \end{align}  
 Using  \eqref{C3} and \eqref{C6}, we rewrite \eqref{C2} in the form
 \begin{align}\label{C8}&
 u(\mu) \breve F(\mu)=-\I\big(I_{m_3}+\b(\mu)\big)\begin{bmatrix} \cll_2^* & 0 & \cll_1^* & 0\end{bmatrix},
\\ \label{C9}& 
 \breve F(\mu):=\begin{bmatrix} \mathfrak{A}_2(\mu_2) &\breve F_{12}
\\
\breve  F_{21} & \mathfrak{A}_1(\mu_1)
 \end{bmatrix},
\\ \label{C10}&
\breve F_{21}=-\I\left(\wh \G_2\Pi_1-\begin{bmatrix} 0 & 0
\\
K_{11} & 0
 \end{bmatrix}\right), \quad
 \breve F_{12}=-\I\left(\wh \G_1\Pi_2-\begin{bmatrix} 0 & 0
\\
K_{12} & 0
 \end{bmatrix}\right).
  \end{align}  
According to \eqref{E}--\eqref{E3} and to the equalities $\wh \G_p \Pi_k=\wh \Pi_p \G_k= \wh \Pi_pR\Pi_k$,  we have
 \begin{align}\label{C11}&
 E(\mu)=\breve F(\mu)-\I \breve \cll \breve \cll^*, \quad \breve \cll:=\col \begin{bmatrix} 0 & \cll_2
& 0 & -\cll_1
 \end{bmatrix}.
  \end{align}   
 Taking adjoints to the expressions in \eqref{D34}, we derive $\G \breve \cll=0$. Hence, relations \eqref{D28}, \eqref{D27}
 and \eqref{D28++} yield $u(\mu)\breve \cll \equiv 0$. Thus, taking into account \eqref{C11}, we rewrite \eqref{C8}
 in the form
 \begin{align}\label{C12}&
 u(\mu) E(\mu)=-\I\big(I_{m_3}+\b(\mu)\big)\begin{bmatrix} \cll_2^* & 0 & \cll_1^* & 0\end{bmatrix},
 \end{align}   
and the theorem's statement follows.  
\end{proof}
\begin{Rk} Proposition \ref{simprep} and Theorems \ref{MbTm} and \ref{bM2Tm}  enable us to construct
$\om(\la,\mu)$ and so the reflection coefficient $\rho(\la,\mu)$ using essentially less ``information"
than before.  In particular, it is easy to see that $E_{12}$ is easily expressed via $G_{12}$,
$E_{21}$ is easily expressed via $G_{21}$, and vice versa.
\end{Rk}

\section{The DSTU and self-adjoint subclasses of  the block TBT-matrices, and  $3$-D Toeplitz matrices}\label{3D}
\setcounter{equation}{0}
\subsection{DSTU matrices}
{\bf 1.} Introduce the matrices
\begin{align}\label{F0:}&
 U_p:=\{\delta_{m_p-i-k+1}\}_{i,k=1}^{m_p} \quad (p=1,2,3), \quad U:=U_1 \otimes U_2 \otimes U_3,
 \end{align} 
 where $\delta_s$ is the Kronecker  delta. 
Recall that by the acronym  DSTU (double structured Toeplitz matrices with the property ``U") we denote the subclass of  the block TBT-matrices
such that the equality
\begin{align}\label{F0+}&
U_3 t_j^{(r)}U_3=(t_j^{(r)})^{\tau}
 \end{align} 
is valid.  In particular, \eqref{F0+} holds for the $3$-D Toeplitz matrices.
 It is easy to see that the TBT-structure and the property \eqref{F0+} also yield the important equalities
 \begin{align}\label{F2}&
\clt_r^{\,\tau}=(U_2\otimes U_3) \clt_r(U_2\otimes U_3), \quad T^{\tau}=UTU.
 \end{align}  
 (Recall  Remark \ref{RkTau} where such relations are discussed.)
 We also have
  \begin{align}\label{F2+}&
M_{2p}^{\tau} = M_{3p}, \quad UM_{3p}(U_k\otimes U_3)=M_{3p},
 \end{align}  
 where $k,p\in \{1,2\}, \,\, k\not= p$.   Now, we express $u_k$ via $\wh u_k$.
 
{\bf 2.} Definitions \eqref{D13}--\eqref{D12} and relations \eqref{F0+}--\eqref{F2+} imply that
\begin{align}\label{F3}&
M_{1p}+UM_{4p}^{\tau}(U_k\otimes U_3)=TM_{3p}, \quad M_{4p}+(U_k\otimes U_3)M_{1p}^{\tau}U=M_{2p}T.
 \end{align}   
 From the second equality in \eqref{F2}, it follows that
\begin{align}\label{F3+}&
R^{\tau}=URU, 
 \end{align}   
where $R=T^{-1}$. Now, we are ready to prove the relations
  \begin{align}\label{F4}&
\wh \G_p^{\tau}=U\G_p \wt U_k+\begin{bmatrix} 0 & M_{3p}\end{bmatrix}, \quad \wt U_k:=\begin{bmatrix} 0 & - U_k\otimes U_3
\\ U_k\otimes U_3 & 0
\end{bmatrix}.
 \end{align}  
Indeed, in view of \eqref{D17} and \eqref{F2+}--\eqref{F3+} we have
\begin{align}\nn
\wh \G_p^{\tau}&=URU\begin{bmatrix} M_{2p}^{\tau} & M_{4p}^{\tau} \end{bmatrix}
=URU\begin{bmatrix} -M_{4p}^{\tau}(U_k\otimes U_3) & M_{2p}^{\tau} (U_k\otimes U_3)\end{bmatrix}\wt U_k
\\ \label{F5} &
=UR\begin{bmatrix} M_{1p}-TM_{3p} & M_{3p} \end{bmatrix}\wt U_k=U\G_p\wt U_k+\begin{bmatrix} 0 & UM_{3p}  (U_k\otimes U_3)\end{bmatrix},
 \end{align}   
and \eqref{F4} follows.

Recall equality \eqref{F0} which yields 
\begin{align}\label{F6}&
-\left(\frac{\mu_p+\frac{\I}{2}}{\mu_p-\frac{\I}{2}}\right)^{m_p}\1_{m_p}^*\big(\cla_p^*-\mu_p I_{m_p}\big)^{-1}=\1_{m_p}^*\big(\cla_p^*+\mu_p I_{m_p}\big)^{-1}U_p.
 \end{align}
Hence, we derive
\begin{align}\nn &
\cll^*\big(A_1^*-\mu_1 I\big)^{-1}\big(A_2^*-\mu_2 I\big)^{-1}
\\ & \nn
=\1_{m_1}^*\big(\cla_1^*-\mu_1 I_{m_1}\big)^{-1}\otimes \1_{m_2}^*\big(\cla_2^*-\mu_2 I_{m_2}\big)^{-1}
\otimes I_{m_3}
\\ \nn &
=q(\mu)
\Big(\1_{m_1}^*\big(\cla_1^*+\mu_1 I_{m_1}\big)^{-1}U_1\otimes \1_{m_2}^*\big(\cla_2^*+\mu_2 I_{m_2}\big)^{-1}U_2\otimes I_{m_3}\Big)
\\ & \label{F7}
=q(\mu)\cll^*\big(A_1^*+\mu_1 I\big)^{-1}\big(A_2^*+\mu_2 I\big)^{-1}\big(U_1\otimes U_2\otimes I_{m_3}\big),
 \end{align}
 where
\begin{align}
 & \label{F8}
q(\mu):= \left(\frac{\mu_1-\frac{\I}{2}}{\mu_1+\frac{\I}{2}}\right)^{m_1}\left(\frac{\mu_2-\frac{\I}{2}}{\mu_2+\frac{\I}{2}}\right)^{m_2}.
 \end{align}
 \begin{Tm}\label{Tmu} Let a block TBT-matrix $T$ satisfy \eqref{F0+} $($i.e., let it belong to the DSTU subclass$)$ and be invertible.
 Then, the matrix functions $u_p$ are expressed via $\wh u_p$ by the formulas
\begin{align}\label{F9}
u_p(\mu)=-q(\mu)U_3\wh u_p(\mu)^{\tau} \wt U_k \quad (k,p\in \{1,2\}, \,\, k\not=p),
\end{align} 
where $u_p$ and $\wh u_p$ are defined by \eqref{D28}--\eqref{D28++}, and $\wt U_k$ are introduced in \eqref{F4}. Moreover, for $G(\la)$ given by \eqref{g}--\eqref{g3}
$($and used in order to construct $\wh u)$  we have
\begin{align}\label{F14}
G_{21}=\wt U_1 G_{12}^{\tau} \wt U_2.
\end{align} 
 \end{Tm}
\begin{proof}. From \eqref{D28++}, \eqref{F4} and \eqref{F7} we derive
\begin{align}\nn
\wh v_p(\mu)^{\tau}=&\cll^*\Big(\prod_{i=1}^2(A_i^*+\mu_i I)^{-1}\Big)\wh \G_p^{\tau}=\cll^*\Big(\prod_{i=1}^2(A_i^*+\mu_i I)^{-1}\Big)U \G_p\wt U_k
\\ \nn &
+\cll^*\Big(\prod_{i=1}^2(A_i^*+\mu_i I)^{-1}\Big)\begin{bmatrix} 0 & M_{3p}\end{bmatrix}
\\ \nn
=&U_3\cll^*\Big(\prod_{i=1}^2(A_i^*-\mu_i I)^{-1}\Big)\G_p\wt U_k\big/ q(\mu_1,\mu_2)+\1_{m_p}^*(\cla_p^*+\mu_pI_{m_p})^{-1}\1_{m_p}
\\ \label{F10} &
\times \big(\1_{m_k}^*(\cla_k^*+\mu_kI_{m_k})^{-1}\otimes I_{m_3}\big)\begin{bmatrix} 0 & I_{m_km_3}\end{bmatrix}.
\end{align} 
In view of \eqref{F0}, we have $\1_{m_p}^*(\cla_p^*+\mu_p)^{-1}\1_{m_p}
=\I\left(1-\Big(\frac{\mu_p+\frac{\I}{2}}{\mu_p-\frac{\I}{2}}\Big)^{m_p}\right).$ Hence, taking into account
\eqref{F6} we rewrite the second term in \eqref{F10} in the form
\begin{align}\nn &
\1_{m_p}^*(\cla_p^*+\mu_pI_{m_p})^{-1}\1_{m_p}
\big(\1_{m_k}^*(\cla_k^*+\mu_kI_{m_k})^{-1}\otimes I_{m_3}\big)\begin{bmatrix} 0 & I_{m_km_3}\end{bmatrix} \\
& \nn
=\I \big(\1_{m_k}^*(\cla_k^*+\mu_kI_{m_k})^{-1}\otimes I_{m_3}\big)\begin{bmatrix} 0 & I_{m_km_3}\end{bmatrix}
\\ & \quad \label{F11} 
+
\I \big(\1_{m_k}^*(\cla_k^*-\mu_kI_{m_k})^{-1}U_k\otimes I_{m_3}\big)\begin{bmatrix} 0 & I_{m_km_3}\end{bmatrix}\big/q(\mu).
\end{align} 
The definitions \eqref{D28+}, \eqref{D27+} and \eqref{D28++} of $\wh u_p$ and  $v_p$  as well as equalities \eqref{F10} and \eqref{F11} yield
\begin{align}\nn
\wh u_p(\mu)^{\tau}=&\frac{1}{q(\mu)}\big(U_3v_p(\mu)\wt U_k
+
\I \big(\1_{m_k}^*(\cla_k^*-\mu_kI_{m_k})^{-1}U_k\otimes I_{m_3}\big)\begin{bmatrix} 0 & I_{m_km_3}\end{bmatrix}\big)
\\ \nn =& 
\frac{U_3}{q(\mu)}\Big(v_p(\mu)\wt U_k
+
\I \big(\1_{m_k}^*(\cla_k^*-\mu_kI_{m_k})^{-1}\otimes I_{m_3}\big)\begin{bmatrix} 0 & U_k\otimes U_3\end{bmatrix}\Big)
\\ \nn =& 
\frac{U_3}{q(\mu)}\Big(v_p(\mu)
-
\I \big(\1_{m_k}^*(\cla_k^*-\mu_kI_{m_k})^{-1}\otimes I_{m_3}\big)\begin{bmatrix} I_{m_km_3} & 0\end{bmatrix}\Big)\wt U_k.
\end{align} 
 Now, \eqref{F9} is immediate from the definitions of $u_p$.

{\bf 3.}  Next, we will prove \eqref{F14}. Using \eqref{D17} and \eqref{F2+}--\eqref{F3+},  we derive
\begin{align}\nn 
\big(\wh \Pi_1 \G_2\big)^{\tau}&=\Pi_2^{\tau}URU\wh \Pi_1^{\tau}=\begin{bmatrix}  M_{12}^{\tau} \\  M_{22}\end{bmatrix}
URU\begin{bmatrix} M_{31} & M_{41}^ {\tau} \end{bmatrix}
\\ \nn &
=\wt U_1\left(\wh \Pi_2R\Pi_1 -\begin{bmatrix} 0 \\ M_{22}T \end{bmatrix}R\Pi_1-\wh \Pi_2 R \begin{bmatrix} TM_{31} & 0 \end{bmatrix}\right.
\\ \label{F15} & \quad 
\left.+\begin{bmatrix} 0 & 0 \\ M_{22} T M_{31} & 0 \end{bmatrix}
\right)\wt U_2
\\ \nn &
=\wt U_1\left(\wh \Pi_2 \G_1 -\begin{bmatrix} M_{22}M_{31}  & 0 \\ M_{22}M_{11}+M_{42}M_{31}  -M_{22} T M_{31} & M_{22}M_{31} \end{bmatrix}
\right)\wt U_2.
\end{align}
According to \eqref{D14} and \eqref{D9}, we have
\begin{align}\nn
M_{22}M_{31}&=\1_{m_1}\otimes \1_{m_2}^*\otimes I_3=(\1_{m_1} \1_{m_2}^*)\otimes I_3=\cll_1\cll_2^*
\\ \label{F16}&
=(U_1\otimes U_3)\cll_{1}\cll_2^*(U_2\otimes U_3)
 \end{align}
Taking into account the first equality in \eqref{F3}, we obtain
\begin{align}\label{F17}&
M_{22}M_{11}+M_{42}M_{31}  -M_{22} T M_{31} =M_{42}M_{31}-M_{22}UM_{41}^{\tau}(U_2\otimes U_3).
 \end{align}
 We also note that
 \begin{align}\label{F18}&
 (U_1\otimes U_3)K_{12}^{\tau}(U_2\otimes U_3)=\begin{bmatrix} \frac{1}{2}\wt M_{12}^{(0)} +\sum_{j=1}^{m_1-1}\wt M_{12}^{(-j)} 
 \\ \frac{1}{2}\wt M_{12}^{(0)} +\sum_{j=1}^{m_1-2}\wt M_{12}^{(-j)} \\ \ldots \\ \frac{1}{2}\wt M_{12}^{(0)}
 \end{bmatrix},
 \end{align}
where $K_{12}$ has the form \eqref{g4} and
 \begin{align}\label{F19}&
\wt M_{12}^{(r)} = \begin{bmatrix}\frac{1}{2} t_{0}^{(r)} +\sum_{j=1}^{m_2-1} t_{j}^{(r)} & \frac{1}{2} t_{0}^{(r)} +\sum_{j=1}^{m_2-2} t_{j}^{(r)} & \ldots & 
\frac{1}{2} t_{0}^{(r)} \end{bmatrix}.
 \end{align}
In view of \eqref{F19}, we have
 \begin{align}\label{F20}&
\wt M_{12}^{(r)}+ M_{42}^{(r)} =\big(\1_{m_2^*}\otimes I_3\big)\clt_r.
 \end{align}
Formulas \eqref{F0+} and \eqref{F18}--\eqref{F20} imply that
 \begin{align}\label{F21}&
(U_1\otimes U_3)K_{12}^{\tau}(U_2\otimes U_3)+M_{42}M_{31}=K_{11}+M_{22}UM_{41}^{\tau}(U_2\otimes U_3).
 \end{align}
By virtue of \eqref{F17} and \eqref{F21} we obtain
\begin{align}\label{F22}&
M_{22}M_{11}+M_{42}M_{31}  -M_{22} T M_{31} +(U_1\otimes U_3)K_{12}^{\tau}(U_2\otimes U_3)=K_{11}.
 \end{align}
Formulas \eqref{g3}, \eqref{F15}, \eqref{F16} and \eqref{F22} yield
\begin{align}\nn
\wt U_1 G_{12}^{\tau} \wt U_2=&\I \left(\wh \Pi_2 \G_1 -\begin{bmatrix} M_{22}M_{31}  & 0 \\ M_{22}M_{11}+M_{42}M_{31}  -M_{22} T M_{31} & M_{22}M_{31} \end{bmatrix}
\right.
\\ \nn &
\left. - \wt U_1 \begin{bmatrix} \cll_1\cll_2^* & K_{12}^{\tau} \\ 0 & 0\end{bmatrix} \wt U_2\right)
\\ \label{F23} =&\I \left(\wh \Pi_2 \G_1 -\begin{bmatrix}\cll_1\cll_2^*  & 0 \\ K_{11} & 0 \end{bmatrix}
\right)=G_{21}.
 \end{align}
\end{proof}
\begin{Rk}\label{RkG} In the process of 
obtaining $\om$ or, equivalently, the reflection coefficient $\rho$ $($for $T\in$ \rm{DSTU}$)$,
we can recover $G$ from  $G_{21}$ instead of $G_{12}$ and we can recover $\wh u$ from $u$
in the same way as we recover $u$ from $\wh u$.
\end{Rk}
\subsection{Self-adjoint block TBT-matrices}
Self-adjoint matrices $T=T^*$ present the most important subclass of the block TBT-matrices.
For this subclass we have
\begin{align}\label{S1}&
T=T^*, \quad R=R^*, \quad \clt_r^*=\clt_{-r}, \quad \big(t^{(r)}\big)^*=t_{-j}^{(-r)},
 \end{align}
where the second equality is immediate from the first one,
the third  equality follows from the first one and from \eqref{D0},
and the last equality is implied by the third one and by \eqref{D0}.
Using \eqref{S1} and definitions \eqref{D13}, \eqref{D15},
we derive
\begin{align}\label{S2}&
\big(M_{11}^{(r)}\big)^*=M_{41}^{(r)}, \quad M_{11}^*=M_{41}.
 \end{align}
 Taking into account \eqref{S1} and definitions \eqref{D9}--\eqref{D12}, we obtain
 \begin{align}\label{S3}&
\big(M_{12}^{(r)}\big)^*=M_{42}^{(-r)}, \quad M_{12}^*=M_{42}.
 \end{align}
Relations \eqref{S1}--\eqref{S3} yield
 \begin{align}\label{S4}&
\big(\wh \Pi_1 R \Pi_2\big)^*=\Pi_2^* R \wh \Pi_1^*= \begin{bmatrix} M_{42} \\ M_{22}\end{bmatrix}R 
\begin{bmatrix} M_{31} & M_{22}\end{bmatrix}=J_2\wh \Pi_2 R \Pi_1 J_1,
 \end{align}
 where 
 \begin{align}\label{S5}&
J_1=\begin{bmatrix} 0 & I_{m_2m_3} \\ I_{m_2m_3} & 0\end{bmatrix}, \quad  J_2=\begin{bmatrix} 0 & I_{m_1m_3} \\ I_{m_1m_3} & 0\end{bmatrix}.
 \end{align} 
 Now, we are ready to formulate and prove a theorem for the  self-adjoint case.
  \begin{Tm}\label{TmSA} Let a block TBT-matrix $T=T^*$  be invertible.
 Then, the matrix functions $u_p$ are expressed via $\wh u_p$ by the formulas
\begin{align}\label{S6}
u_p(\mu)=\wh u_p(\ov{\mu})^*J_p \quad (p=1,2; \,\, \ov{\mu}=(\ov{\mu_1},\, \ov{\mu_2}),
\end{align} 
where $u_p$ and $\wh u_p$ are defined by \eqref{D28}--\eqref{D28++} and the matrices $J_p$ are introduced in \eqref{S5}. Moreover, for $G(\la)$ given by \eqref{g}--\eqref{g3}
$($and used in order to construct $\wh u)$  we have
\begin{align}\label{S7}
G_{21}=-J_2\left(G_{12}^*-\I \begin{bmatrix} \cll_1 \cll_2^* & 0 \\ 0 &  -\cll_1 \cll_2^* \end{bmatrix}\right)J_1.
\end{align} 
 \end{Tm}
 \begin{proof}. The second equalities in \eqref{S2} and \eqref{S3} imply that
\begin{align}\label{S8}
\Pi_p= \wh \Pi_p^*J_p \quad (p=1,2).
\end{align}  
 Hence, using the definitions \eqref{D17} and  \eqref{D28++}, and  the equality $R=R^*$, we easily derive
 \begin{align}\label{S9}
v_p(\mu)=\wh v_p(\ov{\mu})^*J_p,
\end{align} 
and \eqref{S6} follows. 

According to \eqref{g3} and \eqref{S4} we have
 \begin{align}\label{S10}
-J_2G_{12}^*J_1=\I\left(\wh \Pi_2 R\Pi_1-\begin{bmatrix} 0 & 0 \\ K_{12}^* &  \cll_1 \cll_2^* \end{bmatrix}\right).
\end{align} 
Relations \eqref{g5}, \eqref{g4} and \eqref{S3} yield $K_{12}^*=K_{11}$, and so \eqref{S10} may be rewritten in the form
 \begin{align}\label{S11}
G_{21}=-J_2G_{12}^*J_1+\I \begin{bmatrix} -\cll_1 \cll_2^* & 0 \\ 0 &  \cll_1 \cll_2^* \end{bmatrix}.
\end{align}
Formula \eqref{S7} immediately follows from \eqref{S11}.
 \end{proof}

\subsection{$3$-D Toeplitz matrices}
It is known (see, e.g., \cite{SaA73}) that Toeplitz matrices have the property-$U$. Thus, $3$-D Toeplitz matrices
described by the relations \eqref{D1} satisfy \eqref{F0+} and belong to the subclass DSTU.
\begin{Cy}\label{Spc} Let a   $3$-D Toeplitz matrix $T$  be invertible.
 Then, the matrix functions $u_p$ are expressed via $\wh u_p$ by the formulas \eqref{F9}.
 The matrix function $G(\la)$, which is used for the  construction of $\wh u$, is determined
 by $G_{12}$ and \eqref{F14} holds.
\end{Cy}
The third matrix identity for the $3$-D Toeplitz matrices has the form
\begin{align}\label{D4'}&
A_3T-TA_3^*=\I\big(M_{13}M_{23}+ M_{33}M_{43}\big), \quad A_3=I_{m_1m_2}\otimes\cla_3,
 \end{align}
where
\begin{align}\label{D5}&
M_{13}=\{M_{13}^{(i-k)}\}_{i,k=1}^{m_1}, \quad M_{13}^{(r)}=\{M_{13}^{(r,i-k)}\}_{i,k=1}^{m_2} , \quad M_{23}=I_{m_1 m_2} \otimes \1_{m_3}^*, \\
\label{D6}&
M_{33}=M_{23}^*, \quad M_{43}=\{M_{43}^{(i-k)}\}_{i,k=1}^{m_1}, \quad M_{43}^{(r)}=\{M_{43}^{(r,i-k)}\}_{i,k=1}^{m_2},
\\ \label{D7}&
 M_{13}^{(r,s)}=
\col \begin{bmatrix} \tau_0^{(r,s)}/2 &  (\tau_0^{(r,s)}/{2})+\tau_1^{(r,s)} &\ldots & ({\tau_0^{(r,s)}}/{2})+\displaystyle{\sum_{j=1}^{m_3-1}}\tau_{j}^{(r,s)}\end{bmatrix},
\\ \label{D8}&
M_{43}^{(r,s)}=\begin{bmatrix} \tau_0^{(r,s)}/2 & ({\tau_0^{(r,s)}}/{2})+\tau_{-1}^{(r,s)}& \ldots & ({\tau_0^{(r,s)}}/{2})
+\displaystyle{\sum_{\ell=1}^{m_3-1}}\tau^{(r,s)}_{-\ell}\end{bmatrix}.
\end{align}
\begin{Rk} We note that a $3$-D Toeplitz matrix is determined by \\ $(2m_1-1)(2m_2-1)(2m_3-1)\approx 8m_1m_2m_3$ entries.
The entries of $G_{12}$ $($or $G_{21})$, which we need to recover $G$, and the entries of $K$ $($required to recover $\theta$ in \eqref{T1}$)$
present together $5 m_1m_2m_3^2$ ``information" entries. Here, the power $m_3^2$ is caused by the fact that the identity \eqref{D4'} is not used
sufficiently well.
\end{Rk}
{\bf Open problem II.} Is there a way to use three matrix identities \eqref{D4} simultaneously
and thus further reduce the number of ``information" entries, which are necessary in order
to recover $T^{-1}$?

\vspace{0.2em}

{\bf Acknowledgments.}  {The research of Alexander Sakhnovich    was supported by the
Austrian Science Fund (FWF) under Grant  No. P29177.}


\begin{flushleft}

Inna Roitberg, \\
Maschinenbauabteilung TGM,
Vienna,
Austria, \\
e-mail: {\tt innaroitberg@gmail.com}

\vspace{0.2em}

Alexander Sakhnovich, \\
Faculty of Mathematics,
University
of
Vienna, \\
Oskar-Morgenstern-Platz 1, A-1090 Vienna,
Austria, \\
e-mail: {\tt oleksandr.sakhnovych@univie.ac.at}

\end{flushleft}

\end{document}